\documentclass[12pt,a4paper]{article}
\usepackage{amsmath,amsfonts,amssymb,amsthm,cite,mathtools}
\usepackage{graphicx,xcolor}
\usepackage[normalem]{ulem}
\usepackage{times}

\theoremstyle{definition}                    

\theoremstyle{remark}

\newcommand{\C}{\mathbb{C}}

\renewcommand{\H}{\mathcal{H}}      

\newcommand{\N}{\mathbb{N}}         
\newcommand{\R}{\mathbb{R}}         

\renewcommand{\SS}{\mathcal{S}}     

\newcommand{\Th}{\Theta}

\renewcommand{\.}{\cdot}            

\usepackage{xcolor}

\newcommand{\cv}{= \vcentcolon}

\newcounter{mnotecount}[section]
\renewcommand{\themnotecount}{\thesection.\arabic{mnotecount}}
\newcommand{\mnote}[1]
{\protect{\stepcounter{mnotecount}}$^{\mbox{\footnotesize
$
\bullet$\themnotecount}}$ \marginpar{
\raggedright\tiny\em
$\!\!\!\!\!\!\,\bullet$\themnotecount: #1} }

\definecolor{darkgreen}{rgb}{0,.5,0}

\renewcommand{\title}[1]{\vspace{10mm}\noindent{\Large{\bf #1}}\vspace{8mm}}
\newcommand{\authors}[1]{\noindent{\large #1}\vspace{3mm}}
\newcommand{\address}[1]{{\itshape #1\vspace{2mm}}}

\begin{document}

\begin{center}
\title{The Moyal Sphere}

\authors{Micha\l{} Eckstein$^1$, Andrzej Sitarz$^{1,2}$, Raimar Wulkenhaar$^3$}
 
\address{$^{1}$\,Institute of Physics, 
Jagiellonian University, \\
\L{}ojasiewicza 11, 30-348 Krak\'ow, Poland.}

\address{$^{2}$\, Institute of Mathematics 
of the Polish Academy of Sciences, \\
\'Sniadeckich 8, 00-950 Warszawa, Poland.}

\address{$^{3}$\,Mathematisches Institut der Westf\"alischen
  Wilhelms-Universit\"at\\
Einsteinstra\ss{}e 62, D-48149 M\"unster, Germany}
\end{center}

\begin{abstract}
We construct a family of constant curvature metrics on the Moyal plane and
compute the Gauss--Bonnet term for each of them. They arise from the conformal
rescaling of the metric in the orthonormal frame approach. We find a particular 
solution, which corresponds to the Fubini--Study metric and which equips 
the Moyal algebra with the geometry of a noncommutative sphere.
\end{abstract}


\tableofcontents

\section{Introduction}

Noncommutative geometry provides a unified framework to describe
classical, discrete as well as singular or deformed spaces
\cite{ConnesBook,Polaris}.  Most of the examples studied so far were
constructed with a fixed metric, allowing only small pertubations of
the gauge type. Only recently a class of models with conformally
modified noncommutative metrics was constructed for the
noncommutative tori \cite{CoTr} allowing for the computation of the
scalar curvature using different approaches
\cite{CoMo11,FaKh11,DaSi13,Ro13}.

The Moyal deformation of a plane \cite{BFFLS,Rieffel} is one of the
oldest and best-studied noncommutative spaces. Motivated by the
appearence in the string-theory target space \cite{SeWi99} it is often
used as model of noncommutative space-time (see \cite{Sz03} for a
review) and a background for a noncommutative field theory
\cite{Wu06}. However, so far the only geometry considered is the flat
geometry, which corresponds to the constant metric on the plane. In
this paper we introduce a class of conformally rescaled metrics on the
two-dimensional Moyal plane using the orthonormal frame formalism
adapted to the noncommutative setting. We compute the scalar curvature
and look for the solutions of the constant curvature condition, thus
finding the Moyal--Fubini--Study metrics.

The paper is organized as follows: first, we recall the classical (commutative)
constant curvature solutions, then we briefly review the noncommutative
Moyal plane and compute the scalar curvature using conformally 
rescaled orthonormal frames. We discuss explicit solutions in the matrix basis for
the Moyal algebra as well the first order perturbative correction to the 
Fubini--Study metric using smooth functions on the plane. 

\section{Classical Fubini--Study metric on the plane}

Consider the conformally rescaled metric on the plane:

$$ k(x,y)^2 (dx^2 + dy^2), $$

We assume that $k=k(r)$, then the scalar curvature is:

$$ R(k) = 2 k(r)^{-4} \left( k'(r) k'(r) - \frac{k(r)}{r} k'(r)  - k(r) k''(r)  \right) $$

Now, let us look for $k$ such that $R(k) = C =\hbox{const}$.
We obtain the differential equation:

$$ k(r) k''(r) + \frac{k(r)}{r} k'(r) - k'(r) k'(r) = C k^4(r), $$

which has a family of nondegenerate solutions for $A>0$, $a \geq 1$ and $b>0$:

$$ k(r) = {A r^{a-1} \over b + r^{2a} }, $$

so that the scalar curvature is:

$$ R(k) =  {8 a^2 b \over A^2}, $$

the volume:

$$ V(k) = \pi {A^2 \over ba } $$

and the Gauss--Bonnet term:

$$ \int \sqrt{g} R = 8 \pi a. $$

More generally, assume that we have a class of conformally rescaled metrics for which 
the asymptotics of $k(r)$ at infinity is given by:
$$ k(r) \sim r^{-\alpha} + C_1 r^{-\alpha-1} + C_2 r^{-\alpha-2}$$

One can easily verify that the scalar curvature is regular, that is: 
$$ \lim_{r \to \infty} R(r) < \infty, $$
if $C_1 = 0$ and $\alpha \leq 2$. 

We can now compute the Gauss--Bonnet term for such metrics:
$$
\begin{aligned}
 & \int_{\R^2} \sqrt{g} R(g) = 
4\pi \int_0^\infty r dr \; k(r)^{-2} \left( k'(r) k'(r) - \frac{k(r)}{r} k'(r)  - k(r) k''(r)  \right) \\
& = 4\pi \int_0^\infty \!\!\!\! dr \, \left(  r k(r)^{-2} k'(r) k'(r) -  k(r)^{-1} k'(r)  -  r k(r)^{-1} k''(r)  \right) \\
& = 4\pi \int_0^\infty \!\!\!\! dr \, ( - r k'(r) k(r)^{-1} )' = 4 \pi \left( ( - r k'(r) k(r)^{-1} )_\infty - ( - r k'(r) k(r)^{-1} )_0 \right).
\end{aligned}
$$
Assuming that $k(r)$ and its derivatives are regular at $r=0$ and that $k(r)$ behaves like 
$r^{-\alpha}$ at $r \to \infty$ we obtain:
$$ \int_{\R^2} \sqrt{g} R(g) = 4 \pi \alpha. $$
Note that the special case, where $\alpha=2$, which is the Fubini--Study metric yields
the correct Gauss--Bonnet term for the sphere.

On the other hand, if we require that the metric alone remains 
bounded at $r=\infty$, we have, after change of variables $\rho=\frac{1}{r}$, that
$$  \lim_{\rho \to 0} k^2 \bigl( \rho^{-1} \bigr) \rho^{-4} < \infty. $$
This happens if the asymptotics of $k(r)$ is $r^{-\alpha}$ for $\alpha \geq 2$. 
Moreover, if we require that the metric does not vanish at $r=\infty$ 
then $\alpha=2$ is the only solution for the asymptotic behavior of $k(r)$.

\section{The flat geometry of the Moyal plane}

We begin by reviewing here shortly basic results on Moyal geometry.  We take the algebra of the
Moyal plane, $A_\theta$, as a vector space $(\SS(\R^{2}), *)$, ($\SS$ is the Schwartz space),
equipped with the Moyal product $*$ defined as follows
through the oscillatory integrals \cite{Rieffel}:
\begin{eqnarray}
\label{mogeneral}
(f * g)(x):=(2\pi)^{-2}\int_{\R^{2}\times\R^{2}}  e^{i\xi(x-y)}
f(x-\textstyle{\frac{1}{2}}\Theta \xi) \;g(y)\;d^ny\; d^n\xi.
\end{eqnarray}
where $ \Theta = \left( \begin{array}{cc} 0 & \theta \\ - \theta & 0 \end{array} \right)$,
$\theta \in \R$. With the product defined in (\ref{mogeneral}) it is easy to see that 
$$ f \mapsto \int_{\R^n} f(x)\, d^nx $$ is a trace on the Moyal algebra and the standard partial 
derivations $\partial_{x_1},\partial_{x_2}$ remain derivations on the deformed algebra.

The algebra can be faithfully represented on the Hilbert space of 
$L^2$-sections of the usual spinor bundle over $\R^2$ (which is $\H=L^2(\R^{2})\otimes\C^{2}$)
acting by Moyal left multiplication $L^\Th(a)$. It was demonstrated first in \cite{Himalia} that 
the Moyal plane algebra with one of its preferred unitizations yield nonunital, real spectral triples 
for the standard Dirac operator on the plane arising from the flat Euclidean metric. 

\subsection{Matrix basis for the Moyal algebra}

It will be convenient to work with the matrix basis for the Moyal algebra \cite{Himalia}. 
We define first: 
$$ f_{0,0}  = 2 e^{- \frac{1}{\theta} (x_2^2 + x_2^2)}, $$
and the algebra $A_\theta$ has a natural basis consisting of:
$$ f_{m,n} = {1 \over \sqrt{ m! \,n! \, \theta^{m+n}}} \, (a^*)^m \ast f_{00} \ast (a)^n, $$

or more explicitly:
\begin{equation}
f_{m,n}(r, \phi) = 2 (-1)^m  \sqrt{\frac{m!}{n!}} e^{i \phi(m-n)}
\left(\sqrt{\frac{2}{\theta}} r\right)^{(n-m)} L^{n-m}_m\left(
  \frac{2r^2}{\theta}\right) e^{-\frac{r^2}{\theta}}. 
\label{Moyal-radial}
\end{equation}

We have for each $m,n,k,l \geq 0$:

$$ f_{m,n} \ast f_{k,l} = \delta_{kn} f_{m,l}, \;\;\; f_{m,n}^* = f_{n,m}.$$

In particular for all $k \geq 0$ all $f_{k,k}$ are pairwise orthogonal
projections of rank one.

The natural (not normalized) trace on the Moyal algebra is:
$$ \tau(f) = \int_{\R^2} d^2x \, f(x), \;\;\; \tau(f_{m,n}) = 2 \pi \theta \delta_{mn}. $$

It is convenient to work with the linear combinations of derivations:

$$ \partial = \frac{1}{\sqrt{2}} (\partial_{x_1} - i \partial_{x_2}), \;\;\;\;
 \bar{\partial} =  \frac{1}{\sqrt{2}} (\partial_{x_1} + i \partial_{x_2}). 
 $$

Then:

$$ 
\partial f_{m,n} = \sqrt{\frac{n}{\theta}} f_{m,n-1} - 
 \sqrt{\frac{m+1}{\theta}} f_{m+1,n}, $$
$$\bar{\partial} f_{m,n} = \sqrt{\frac{m}{\theta}} f_{m-1,n} - 
 \sqrt{\frac{n+1}{\theta}} f_{m,n+1}, $$
 
 Furthermore,
 $$ \partial \bar{\partial} = \frac{1}{2} (\partial_{x_1}^2 + \partial_{x_2}^2 ). $$
 
\subsection{Radial functions in the matrix basis}

We define {\em radial} functions on the Moyal plane as those which in their
presentation in the matrix basis have only diagonal elements $f_{n,n}$,
so:

$$ h = \sum_{n=0}^\infty h_n f_{n,n}. $$

Each function $F$ applied to $h$ is easily computable, 
since $f_{n,n}$ are projections, we have:
$$ F(h) = \sum_{n=0}^\infty F(h_n) f_{n,n}. $$

From the orthogonality of matrix basis $f_{n,n}$, the explicit 
representation (\ref{Moyal-radial}) and known integrals of Laguerre
polynomials one deduces
\[
r^a= \sum_{n=0}^\infty \theta^{\frac{a}{2}}\;
{}_2F_1(-n,-\tfrac{a}{2};1;2)\, f_{n,n}(r).
\]
In particular,
\[
\sum_{n=0}^\infty f_{n,n}=1,\qquad
\sum_{n=0}^\infty n f_{n,n}=\frac{1}{2}\Big(\frac{r^2}{\theta}-1\Big).
\]

Therefore the inverse of the conformal factor for the Fubini--Study metric, which 
is a linear function of $r^2$, $k(r)^{-1} = \frac{1}{A} (b + r^2)$ has the following 
expansion:
$$ \frac{1}{A} (b + r^2) = \frac{1}{A} \sum_{n=0}^\infty \left(2 \theta n + b + \theta  \right) f_{n,n}. $$


\section{Conformally rescaled metric on the Moyal plane}

There are several approaches to the conformal rescaling of the metric and 
computing the curvature for the noncommutative torus and -- by analogy -- 
the Moyal plane. Note that each of them uses a different notation and they 
do not give compatible results in the case of the noncommutative torus. 

First of all, one can take a conformal rescaling of the flat Laplace
operator: $ \Delta_h = h \Delta h$, which in the case of the
noncommutative torus was studied by Connes--Tretkoff \cite{CoTr}. This
corresponds to the rescaling of the metric by $h^{-2}$

A second approach uses the language of orthonormal frames
\cite{DaSi13}. It replaces the standard (flat) orthonormal frames,
understood as derivations on the algebra, by the conformally
rescaled ones, $\delta_i \to h \delta_i$, where $h$ is from the
algebra (or, more generally, from its multiplier). Applying the
classical formula for the scalar curvature, which easily adapts to the
noncommutative case, one obtains a noncommutative version of the scalar
curvature.  This corresponds also to the rescaling of the metric by
$h^{-2}$.

Finally, extending the computations of Jonathan Rosenberg for the
Levi-Civita connection \cite{Ro13} on the noncommutative torus one
might obtain a similar generalization of the expression for the scalar
curvature.

We shall concentrate on the case of orthonormal frames and try to
determine whether there exists a radial conformal rescaling for which the
scalar of curvature is constant.

\subsection{The orthormal frame approach}

Let the orthonormal basis of frames be $e_i = h \delta_i$ for $h$ in
the Moyal algebra or its multiplier.\footnote{Note that the rescaled
  frames are no longer derivations, however, if $h$ is taken from the
  commutant of the algebra (or its multiplier), $e_i$ will be
  derivations from the Moyal algebra into the algebra of bounded
  operators on the Hilbert space \cite{DaSi13}.  Since this does not
  change anything in the computations, for the sake of simplicity we
  work with $h$ from the algebra itself.} We assume that $h$ is
positive and invertible, by $h_\ast^{-1}$ we denote the inverse with
respect to the Moyal product, similarly all powers are also Moyal powers.

We have:
$$ [e_1, e_2] =  h \ast \delta_1(h) \ast h_\ast^{-1} e_2 -   h \ast \delta_2(h) \ast  h_\ast^{-1} e_1. $$
so that 
$$ c_{122} =  h \ast \delta_1(h) \ast h_\ast^{-1},\;\;\;\;  c_{121} =  - h \ast \delta_2(h)\ast h_\ast^{-1}. $$
and
$$ c_{212} =  - h \ast \delta_1(h) \ast h_\ast^{-1},\;\;\;\;  c_{211} =  h \ast \delta_2(h) \ast h_\ast^{-1}. $$
We compute the scalar curvature as in \cite[(2.1)]{DaSi13}:
$$ R = 2 h \ast \delta_i (h \ast \delta_i(h) \ast h_\ast^{-1}) - (h \ast \delta_i(h) \ast h_\ast^{-1})^2 
- \frac{1}{2} (h \ast \delta_i(h) \ast h_\ast^{-1})^2 - \frac{1}{2} (h \ast \delta_i(h) \ast h_\ast^{-1})^2 $$

which after simplifications yields:
\begin{align}
\label{R_frame}
 R &= 
 2 h^2 \ast (\delta_{ii} h) \ast h_\ast^{-1} + 2 h \ast \delta_i(h) \ast \delta_i(h) \ast h_\ast^{-1} \notag \\
& \phantom{=} - 2 h^2 \ast \delta_i(h) \ast h_\ast^{-1} \ast \delta_i(h) \ast h_\ast^{-1} 
- 2 h \delta_i(h) \ast \delta_i(h) \ast h_\ast^{-1} \\
 &= 2 h^2 \ast (\delta_{ii} h) \ast h_\ast^{-1} - 2 h^2 \ast \delta_i(h) \ast h_\ast^{-1} \ast \delta_i(h) \ast h_\ast^{-1}.  \notag 
\end{align}

So, to solve the equation $R(h) = C = \hbox{const}$ we need to solve:
$$
\label{eq_frame}
(\Delta h) - \delta_i(h) \ast h_\ast^{-1} \ast \delta_i(h) = C h_\ast^{-1}, 
$$
where $\Delta$ is the standard flat Laplace operator, $\Delta = \delta_1^2 + \delta_2^2$.

We shall present the proof of the existence of the solution in the matrix basis as well
as compute explicitly the first term of the perturbative expansion for the Fubini--Study metric.

\subsection{Solution in the matrix basis}\label{sol}

{\bf Ansatz:} First we look for a radial solution:
$$ h = \sum_{n=0}^\infty \phi_n f_{n,n}. $$

Using the action of partial derivatives and the Laplace operator on the basis, 

$$ \Delta f_{m,n} = \frac{2}{\theta} 
\left( -(m\!+\!n+\!1) f_{m,n} \!+\! \sqrt{(m\!+\!1)(n\!+\!1)} f_{m\!+\!1,n\!+\!1} \!+\! \sqrt{mn} f_{m\!-\!1,n\!-\!1}  \right),  $$
we obtain the following equation:
$$ 
\begin{aligned}
&   \sum_n \left( -(2n+1) \phi_n - n \phi_n^2 \phi_{n-1}^{-1}  - (n+1) \phi_n^2 \phi_{n+1}^{-1} \right) f_{n,n} \\
& + \sum_n (n+1) (\phi_{n+1} + \phi_n) f_{n+1,n+1} 
   + \sum_n n (\phi_{n-1} + \phi_n) f_{n-1,n-1} \\
& = R \. \phi_n^{-1} f_{n,n},
\end{aligned}
$$
where we have set $R = C \theta$. It yields the following recurrence relation:
\begin{align}\label{recurrence}
 \frac{n+1}{\phi_{n+1}} (\phi_{n+1}^2 - \phi_n^2)  
 + \frac{n}{\phi_{n-1}} (\phi_{n-1}^2  - \phi_n^2)  =  R\. \phi_{n}^{-1},
\end{align}
for $n \geq 1$. Note that although \eqref{recurrence} is of the second order, it has only one degree of 
freedom, since for $n=0$ we have
\begin{align}\label{rec0}
\phi_1^2 - \phi_0^2  = R \frac{\phi_1}{\phi_0}.
\end{align}

The recurrence relation is a quadratic one, solving it for $x=\phi_{n+1}$ we have:
 $$  (n+1) x^2 + x \left( \frac{n}{\phi_{n-1}} (\phi_{n-1}^2  - \phi_n^2) - \frac{R}{\phi_n} \right) 
 - (n+1) \phi_n^2 = 0, $$ and it is easy to see that it has only one
 positive root. It can also be easily seen that the sum of roots is
 positive and since their products is $-\phi_n^2$ then the positive
 root must be bigger than $\phi_n$. Hence the solution will be an
 increasing positive sequence. Since $h$ needs to be a positive
 operator, we should start with an initial value $\phi_0 > 0$ and take
 the positive root at each step to have $\phi_n > 0$ for every $n \in
 \N$. Note that $\phi_0 = 0$ is not allowed in \eqref{recurrence}.

\subsubsection{Asymptotic behavior and Gauss--Bonnet term}

As we have shown, there exists a family of solutions yielding a
positive constant scalar curvature for the Moyal plane. We shall now look
for some special solutions and their asymptotic behavior.  First of
all, observe that in the orthonormal frame formalism we have $\sqrt{g}
= h^{-2}$, so we can take as the noncommutative volume element
$h_\ast^{-2}$.  In order to have a finite
volume, the growth of the sequence $\phi_n$ must be faster than
$\sqrt{n}$ so that $\phi_n^{-2}$ gives a summable series.
 
For each solution of the recurrence relation \eqref{recurrence} we
shall compute now the Gauss--Bonnet term:
$$ \tau( \sqrt{g} \ast R) = \int_{\R^2} \sqrt{g} R, $$
where, again, $\sqrt{g}$ is the noncommutative Moyal volume
element. Note that the expression has no ambiguity because of the
trace property of $\tau$ and we need to compute:

\begin{align*}
\int_{\R^2} \sqrt{g} R &= 
\tau \left( 2 \left( (\Delta h) - \delta_i(h) \ast h_\ast ^{-1} 
\ast \delta_i(h) \right) \ast h^{-1} \right)
\\
& = \frac{2 \pi}{\theta} \int_0^\infty r dr \, \sum_{n = 0}^{\infty} \big\{ \left[ 
-(2n+1) - n \phi_n \phi_{n-1}^{-1}  - (n+1) \phi_n \phi_{n+1}^{-1} \right] f_{n,n} \\
& \qquad + (n+1) (\phi_{n+1} \phi_n^{-1} + 1) f_{n+1,n+1} 
   + n (\phi_{n-1} \phi_n^{-1} + 1) f_{n-1,n-1}  \big\}.
\end{align*}

Now recall that $ \int_0^\infty r dr f_{n,n}(r) = \theta$ for all $n
\in \N$. However, to compute the integral with need to introduce a
cut-off in the series. We thus have
\begin{multline*}
\sum_{n = 0}^{N} \left[ -(2n+1) - n \phi_n \phi_{n-1}^{-1}  - (n+1) \phi_n \phi_{n+1}^{-1} \right] + \sum_{n = 0}^{N+1}  n (\phi_{n} \phi_{n-1}^{-1} + 1) + \\ + \sum_{n = 0}^{N-1} (n+1) (\phi_{n} \phi_{n+1}^{-1} + 1) = (N+1) ( \phi_{N+1} \phi_N^{-1} - \phi_{N} \phi_{N+1}^{-1} ).
\end{multline*}

The expression above has a finite and nonvanishing limit as $N \to
\infty$ if the sequence $h_{N} = \phi_{N} \phi_{N-1}^{-1}$ has a limit
$1$ and $N (h_N - 1)$ has a finite limit. This requirement alone is
not sufficient to determine the asymptotic form of the solution. As
the recurrence relation is highly nonlinear we can only check that
some asymptotics are compatible with the relations as well as the
above requirement for the Gauss--Bonnet term. In particular, possible
asymptotics include the power-growing sequences $\phi_n \sim A n^a$ as
well as power-growing sequences modified by logarithms.

\begin{center}
\begin{figure}
\includegraphics[width=0.5\textwidth]{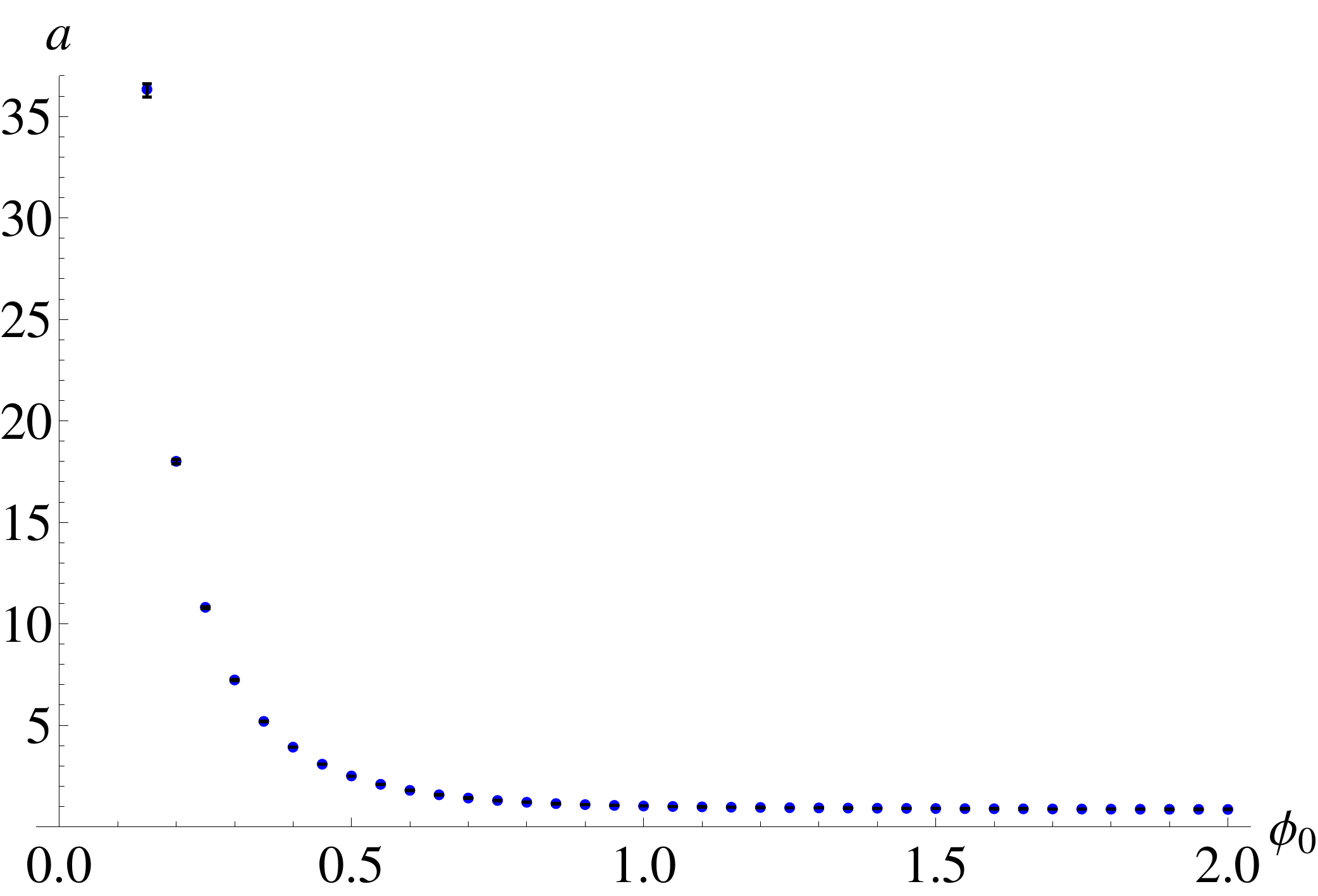} 
\includegraphics[width=0.5\textwidth]{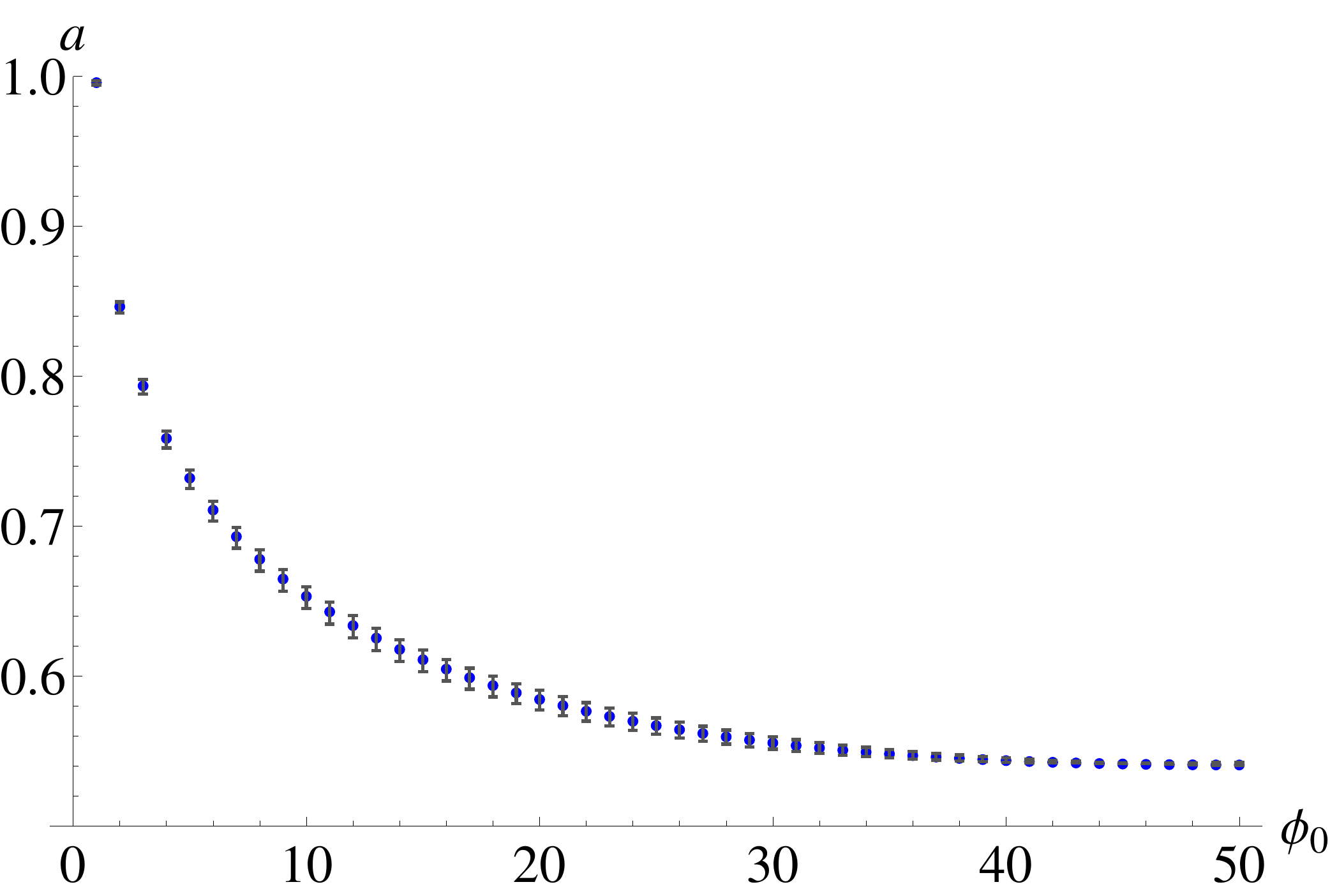}
\caption{\em The above plots were obtained by a numerical computation of $\phi_n$ up to $n = 5000$ starting from 
a given value of $\phi_0$. The asymptotic exponent $a$ is then computed as $\log \phi_N / \log N$ at $N=5000$. The error 
bars are obtained by performing the numerical computation for $N=6000$ (upper) and $N=4000$ (lower). They show the stability 
of the numerical algorithm.\label{fig1}}
\end{figure}
\end{center}

To have an insight into the entire family of possible solutions we
have carried out a numerical study of the solutions (see fig.\ \ref{fig1}), which confirms
that the asymptotics is of that type and gives a glimpse of the
relation $a=a(\phi_0)$.

Assuming that asymptotically $\phi_n \sim A n^a$ as $n$ tends to
$\infty$ for some $A, a > 0$, then
\begin{align*}
\lim_{N \to \infty} (N+1) ( \phi_{N+1} \phi_N^{-1} - \phi_{N} \phi_{N+1}^{-1} ) = 2 a,
\end{align*}
and
\begin{align*}
\int_{\R^2} \sqrt{g} R = 8 \pi a,
\end{align*}
which is, in fact, quite similar to the classical case.

\subsection{The Moyal--Fubini--Study metric}

It is obvious from the above computations that in the case of linear
asymptotics $\phi_n \sim n$, that is $a=1$, we obtain the same Gauss--Bonnet term as for the classical sphere and therefore the solution to \eqref{recurrence} with $\phi_n \sim n$ could be understood
as the Moyal--Fubini--Study solution.

Although we cannot solve exactly the recurrence relation even in this
particular case, one can systematically find the asymptotics of
$\phi_n$. Explicit computations give up to $o(\frac{1}{n^3})$:
$$ \phi_n = n + \frac{1}{2}(R+1) + \frac{1}{8} \frac{1}{n} - \frac{13 R+9}{144}  \frac{1}{n^2}
+ \frac{1}{32} \left( - \frac{1}{4} + \frac{26}{9} R + \frac{29}{18} R^2 \right)  \frac{1}{n^3} + \cdots. $$ 
We remark as well that the case $a=1$ corresponds exactly to the asymptotic behavior 
of the coefficients of the classical Fubini--Study metric. Therefore the Moyal--Fubini--Study is, in 
fact, a perturbation of the classical Fubini--Study metric.

\subsection{Curvature \`a la Rosenberg}

In this section we would like to compare the above computation to the curvature \`a la Rosenberg 
\cite{Ro13}. To recall briefly, in two dimensions we have (classically)
\begin{align*}
R = g^{\mu\nu} R_{\mu \nu} = g^{\mu\nu} g^{\kappa\lambda} R_{\kappa\mu\lambda\nu} = 2 g^{11} g^{22} R_{1212}. 
\end{align*}
In Rosenberg's convention $g^{11} = g^{22} = e^{-h} \cv H^{-1}$ and (see \cite[(4.3)]{Ro13})
\begin{align*}
R_{1212} = -\tfrac{1}{2} \left( \Delta (H) - \delta_i(H) H^{-1} \delta_i(H) \right).
\end{align*}
Note, that since $[H,\delta(H)] \neq 0$, the scalar curvature is not uniquely defined in this framework as we can
always choose $g^{11} R_{1212} g^{22} \neq g^{11} g^{22} R_{1212}$. However, if we are interested in the case
of the constant curvature it does not affect the equation: 
\begin{align}\label{eq_rosenberg}
\Delta (H) - \delta_i(H) H^{-1} \delta_i(H) = -C H^{-2}.
\end{align}
If one takes a radial ansatz for $H$ as we did in Section \ref{sol}, one discovers that \eqref{eq_rosenberg} 
yields a recursion relation very similar to \eqref{recurrence},
\begin{align}\label{R_recurrence}
 \frac{n+1}{\phi_{n+1}} (\phi_{n+1}^2 - \phi_n^2)  
 + \frac{n}{\phi_{n-1}} (\phi_{n-1}^2  - \phi_n^2)  =  - c \phi_{n}^{-2}.
\end{align}

Let us note that both formulae \eqref{R_frame} and \eqref{eq_rosenberg} have correct classical 
(i.e. $[H,\delta(H)] = 0$) limits since
\begin{align*}
& R(g_{ij} = H \delta_{ij})  = H^{-3} \left( \delta_{i}(H) \delta_{i}(H) - H \Delta H  \right) \\
& R(g_{ij} = h^{-2} \delta_{ij})  = 2 \left( h \Delta h - \delta_{i}(h) \delta_{i}(h) \right).
\end{align*}

This suggests that the curvature \`a la Rosenberg is the same as the
one obtained in the orthonormal frame formalism with $H =
h^{-2}$. However, although this is true in the commutative limit, the
formulae do not match exactly in the non-commutative
framework. Indeed,
$$
\begin{aligned}
R_R(H = h^{-2}) &= h (\Delta h) + (\Delta h) h - h \delta_i(h) h^{-1} \delta_i(h)  \\
&\qquad - \delta_i(h) h^{-1} \delta_i(h) h + \delta_i(h) \delta_i(h) - h \delta_i(h) h^{-2} \delta_i(h) h.
\end{aligned}
$$

\section{Perturbative solution for the Moyal--Fubini--Study metric}

In the physics literature the deformation parameter $\theta$ is frequently treated as 
a small perturbation parameter and all possible physical quantities and fields are 
computed up to certain order in $\theta$. We can follow this principle and attempt to
solve the equation \eqref{eq_frame} up to the next leading order around the classical
Fubini--Study solution. We begin with some technical computations.

\subsection{Perturbative expansion of the Moyal product}

Using the formal expression for the Moyal product of two functions on $\R^2$:
$$ (f \ast g)(p) = \left( e^{ \frac{i\theta}{2} \partial_p \partial_q} f(p) g(q) \right)_{|_{p=q}}, $$
we obtain the following perturbative expansion of the Moyal product of two radial functions:

$$ f \ast g =  f\, g  - \frac{\theta^2}{8r} \left( f'' g' + f' g'' \right), $$
where derivatives are with respect to the argument $r$.

As a consequence the following gives the perturbative formula 
for the Moyal inverse of the radial function:

$$ f^{-1}_\ast = f^{-1} + \frac{\theta^2}{4r} \left( (f')^3 f^{-4} - f'' f' f^{-3} \right). $$

We shall also need the formula for the perturbative expansion of the
following product:

$$ \delta_i(f) \ast f^{-1}_\ast \ast \delta_i(f), $$

Explicit computations give:$$ 
\begin{aligned}
& \delta_i(f) \ast f^{-1}_\ast \ast \delta_i(f) =  (f')^2 f^{-1}  \\
& \quad \quad + \theta^2 \left(
    \frac{1}{4} \frac{1}{r^4} (f')^2 \, f^{-1} 
+  \frac{1}{2}  \frac{1}{r^3} (f')^3 \, f^{-2} 
  -  \frac{1}{2} \frac{1}{r^3} f' \, f'' \, f^{-1} 
 +  \frac{1}{r^2} (f')^4 f^{-3}  
   \right. \\
&    \quad \quad \quad
  -  \frac{3}{2}  \frac{1}{r^2} (f')^2 \, f'' \, f^{-2} 
   + \frac{1}{4} \frac{1}{r} (f')^5 f^{-4}
  +  \frac{1}{4}  \frac{1}{r^2} (f') \, f''' \, f^{-1} 
  +  \frac{1}{4} \frac{1}{r^2} (f'')^2 \, f^{-1} \\
&   \left. \quad \quad \quad   
  -  \frac{3}{4}  \frac{1}{r} (f')^3 \, f'' \, f^{-3} 
 +  \frac{1}{4} \frac{1}{r} (f')^2 \, f''' \, f^{-2} 
   +  \frac{1}{2}  \frac{1}{r} f' \, (f'')^2 \, f^{-2} 
   -  \frac{1}{4}  \frac{1}{r} f'' \, f''' \, f^{-1} \right).
\end{aligned}
$$

\subsection{The Moyal--Fubini--Study metric up to order $\theta^2$}

We look for the solution of the constant curvature equation in the form:
$$ h(r) = h_{FS}(r) + \theta^2 \epsilon(r) + o(\theta^2).$$
The classical Fubini--Study solution is $h_{FS}(r) = \eta (1 + r^2) $, using
the pertubative expansion derived above one obtains the following
differential equation in order $\theta^2$:
$$ 
\begin{aligned}
& \epsilon''(r) + \left( \frac{1}{r} - \frac{4r}{1 + r^2} \right) \epsilon'(r) 
+ \frac{4}{(1 + r^2)} \epsilon(r)  + 8 \eta \frac{(1- r^2)}{(1+r^2)^3}  = 0.
\end{aligned}
$$
The most general solution is of the form:
$$ 
\begin{aligned}
\epsilon(r) =&\, C_1  (r^2 - 1) + C_2 ( (r^2 - 1) \log(r) - 2 ) \\
& - \frac{\eta}{3} \frac{1}{1+r^2} 
\biggl( (1-r^4) \log(1+r^2) + 2 (r^4-1) \log(r) + r^2 -2 \biggr),
\end{aligned}§
$$

which, however, might have a logarithmic singularity at $r\!=\!0$. However, setting
one of the integration constants, $C_2 = \frac{2}{3} \eta $, fixes the problem. The second integration constant determines the value of $\epsilon$ at $r\!=\!0$ as $\epsilon(0) = - C_1 - \frac{2}{3} \eta$.

\begin{center}
\begin{figure}
\includegraphics[width=0.9\textwidth]{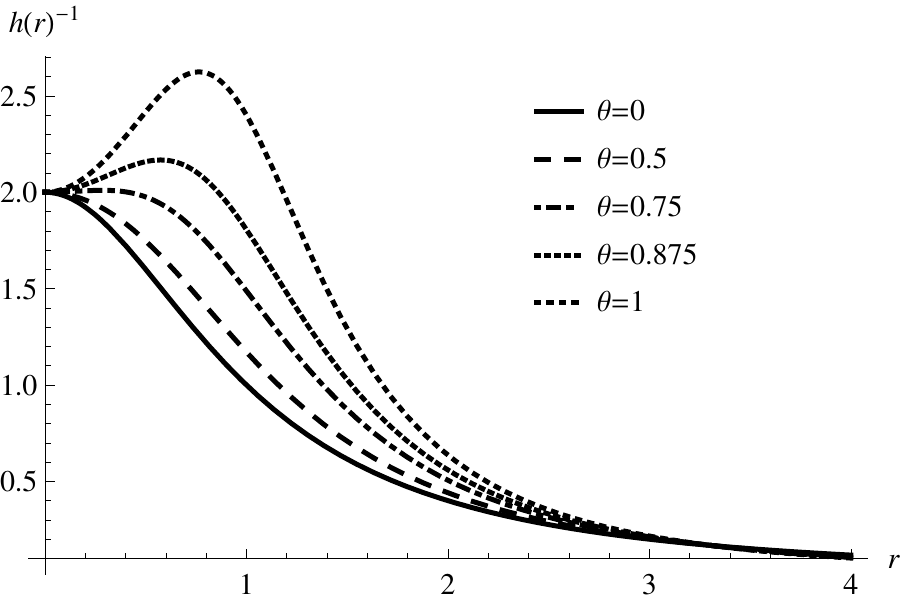}
\caption{\em The graph of the Moyal--Fubini--Study conformal factor $h(r)^{-1}$ up to the order $o(\theta^2)$, for different values of $\theta$. The parameters are chosen $C_1 = -\frac{2}{3}\eta$ (which yields $\epsilon(0)=0$) and $\eta = \frac{1}{2}$ (which corresponds to $R=1$).\label{fig2}}
\end{figure}
\end{center}

In fig.\ \ref{fig2} we present the plot showing how the Fubini--Study conformal
factor gets altered by the deformation parameter $\theta$. Note that
we have displayed the graph of $1/h$, which is equal to the conformal
factor $k$ from Section 2, in the case when $\theta = 0$.

\section{Conclusions and discussion}

We have demonstrated that on the Moyal plane there exists a family of
metrics which resemble the classical Fubini--Study metric and for which 
the scalar curvature computed in the orthonormal frame 
formalism is constant. 

Unlike in the situation of the noncommutative torus some explicit
computation and explicit solutions for the conformally rescaled metric
are possible both in the matrix basis formalism and as a
perturbative solution in $\theta$. Since the $C^\ast$-completion of the Moyal
algebra is the algebra of compact operators, it is natural to consider its simplest
unitization as a model for the quantum sphere like the standard
Podle\'s sphere \cite{Podles}. As the $C^\ast$-algebras describe only
topology this alone cannot distinguish between different geometries,
which could be constructed over these noncommutative topological
spaces.  Therefore only the examples of metrics in fact provide
relevant noncommutative geometries.

The metric we have found would be a natural candidate for 
the ``Moyal sphere'', a spherical noncommutative geometry with
a metric that has constant scalar curvature. This is a basis to
further studies of such geometries and their extensions. First of all,
it would be interesting to see whether the constructed metric leads
to a finitely summable spectral triple of dimension $2$ over the 
unitized Moyal algebra (note that the unitization is different from 
the one in \cite{Himalia}) and to verify whether the distance function 
between the states on the Moyal algebra, in particular vector states 
and coherent states (see \cite{Mar12}) is bounded from above. 

A further question concerns the existence and description of 
a three-dimensional noncommutative sphere with an action
of $U(1)$ group so that the ``Moyal sphere'' is a fixed-point
subalgebra (homogeoneous space). This would lead to the
construction and studies of Moyal magnetic monopole solutions.

As Moyal deformation is used as a model for noncommutative
space-time (in $4$ dimensions) one can check whether similar 
Fubini--Study type solutions and geometries exist in higher 
dimensions.

The sphere-like metric could also provide a method for the regularization of the quantum field theory over the Moyal space. The usual approach
has a severe infrared problem \cite{MiSe} as result of unbounded distances. 
The metric regularization might provide an alternative solution to the 
harmonic oscillator approach \cite{GrWu04}. 

{\bf Acknowledgements:} ME and AS would like to thank the Mathematisches Institut 
der Westf\"alischen Wilhelms-Universit\"at in  M\"unster for hospitality, and RW would 
like to thank the Jagiellonian University in Krak\'ow for hospitality. We thank the 
Foundation for Polish Science IPP Programme ``Geometry and Topology in Physical Models'' 
and the SFB 878 ``Groups, Geometry and Actions'' for funding these mutual visits during 
which most of the work has been done.
  
ME acknowledges the support of the Marian Smoluchowski Krak\'ow Research Consortium 
``Matter--Energy--Future'' within the programme KNOW and the support of the Foundation 
for Polish Science IPP Programme ``Geometry and Topology in Physical Models'' co-financed 
by the EU European Regional Development Fund,Operational Program Innovative Economy 
2007-2013.

AS acknowledges support of NCN grant 2012/06/M/ST1/00169.


\begin{thebibliography}{53}
\bibitem{BFFLS} F. Bayen, M. Flato, C. Fronsdal, A. Lichnerowicz and D. Sternheimer,
{\it Deformation theory and quantization. I. Deformations of symplectic structures,
II. Physical applications}, Ann.\ Phys.\ (NY) {\bf 111} (1978) 61--151.
\bibitem{Mar12} E. Cagnache, F. D'Andrea, P. Martinetti, J.-C. Wallet, 
{\it The spectral distance on the Moyal plane}, 
J.Geom.Phys., \textbf{61}, (2011), 1881-1897,
\bibitem{ConnesBook} A.Connes,
{\it Noncommutative Geometry}, Academic Press 1994.
\bibitem{CoMo11} A. Connes, H. Moscovici, 
{\it Modular curvature for noncommutative two-tori}, 
J.\ Amer.\ Math.\ Soc.\ \textbf{27} (2014) 639. 
\bibitem{CoTr} A.~Connes, P.~Tretkoff, 
{\it The Gauss--Bonnet theorem for the noncommutative two torus},
in: \emph{Noncommutative geometry, arithmetic, and related topics}, 
141--158, Johns Hopkins Univ. Press, Baltimore, MD, 2011
\bibitem{DaSi13} L. Dabrowski and A. Sitarz, 
{\it Curved noncommutative torus and Gauss--Bonnet}, 
Journal of Mathematical Physics \textbf{54} (2013) 013518.
\bibitem{FaKh11} F.~Fathizadeh, M.~Khalkhali, 
{\it Scalar curvature for the noncommutative two torus,}  
 J. Noncommut. Geom. \textbf{7} (2013), 1145--1183.
\bibitem{Himalia} V. Gayral, J.M. Gracia-Bond\'{\i}a, B. Iochum, T. Sch\"ucker
and J.~C. V\'arilly, {\it Moyal planes are spectral triples}, 
Comm.\ Math.\ Phys.\ {\bf 246} (2004) 569--623.
\bibitem{Polaris} J. M. Gracia-Bond\'{\i}a, J. C. V\'arilly and H. Figueroa,
\emph{Elements of Noncommutative Geometry}, Birkh\"auser Advanced Texts,
Birkh\"auser, Boston, 2001.
\bibitem{GrWu04} H. Grosse and R. Wulkenhaar, 
{\it Renormalisation of $\phi^4$-theory on noncommutative $\R^4$ in the matrix base},
Commun. Math. Phys., \textbf{256} (2005) 305-374.

\bibitem{MiSe}
S. Minwalla, M. Van Raamsdonk, N. Seiberg. 
{\it Noncommutative perturbative dynamics},
JHEP 0002, 020 (2000)
\bibitem{Podles} P. Podle\'s
{\it Quantum Spheres}, 
Lett. Math. Phys. \textbf{14} (1987) 193--202
\bibitem{Rieffel} M. Rieffel,
{\it Deformation quantization for actions of $\R^d$},
Memoirs AMS, vol 506, AMS, Providence, (1993).
\bibitem{Ro13} 
J. Rosenberg, {\it Levi-Civita's Theorem for Noncommutative Tori}, 
SIGMA \textbf{9} (2013) 071.
\bibitem{SeWi99} N. Seiberg and E. Witten,
{\it String theory and noncommutative geometry}
JHEP 9909:032, (1999)
\bibitem{Sz03} R. Szabo,
{\it Quantum Field Theory on Noncommutative Spaces},
Phys. Rept. \textbf{378} (2003) 207--299.
\bibitem{Wu06} R. Wulkenhaar,
{\it Field theories on deformed spaces}, 
J. Geom. Phys. \textbf{56} (2006) 108--141.
\end{thebibliography}
\end{document}